\documentclass{amsart}%
\usepackage{amssymb}
\usepackage{amsmath}
\usepackage{latexsym}
\usepackage{hyperref}
\usepackage{amsfonts}
\usepackage{graphicx}%
\newtheorem{theorem}{Theorem}

\newcommand{\be}{\begin{equation}}
\newcommand{\ee}{\end{equation}}
\newcommand{\bea}{\begin{eqnarray}}
\newcommand{\eea}{\end{eqnarray}}

\begin{document}
\title{A necessary and sufficient condition for Ricci shrinkers to have positive AVR}
\author{Bennett Chow}
\address{Department of Mathematics, University of California San Diego, La Jolla, CA 92093}
\email{{benchow@math.ucsd.edu}}
\author{Peng Lu}
\address{Department of Mathematics, University of Oregon, Eugene, OR 97403}
\email{{penglu@uoregon.edu}}
\author{Bo Yang}
\address{Department of Mathematics, University of California San Diego, La Jolla, CA 92093}
\email{{b5yang@math.ucsd.edu}}

\begin{abstract}
In this short note we observe that a recent result of C.-W. Chen meshes well
with earlier work of H.-D. Cao and D.-T. Zhou, O. Munteanu, J. Carrillo and L.
Ni, and S.-J. Zhang to give a necessary and sufficient condition for complete
noncompact shrinking gradient Ricci solitons to have positive asymptotic
volume ratio.

\end{abstract}
\maketitle


Let $\left(  \mathcal{M}^{n},g,f\right)  $ denote a complete shrinking
gradient Ricci soliton (\emph{shrinker} for short) with $R_{ij}+\nabla
_{i}\nabla_{j}f-\frac{1}{2}g_{ij}=0$. Throughout we shall assume that $f$ is
the normalized potential function in the sense that $R+|\nabla f|^{2}-f=0$
holds on $\mathcal{M}$.

It was proved by B.-L. Chen \cite{ChenB} that complete ancient solutions to
the Ricci flow, and in particular shrinkers, must have nonnegative scalar
curvature. As a consequence, the potential function $f$ satisfies the
estimate:%
\begin{equation}
0\leq f(x)\leq\left(  \frac{1}{2}r(x)+f{(O)}^{\frac{1}{2}}\right)  ^{2},
\label{Sharp upper bound for f}%
\end{equation}
where $r(x)$ denotes the distance function to a fixed point $O$ in
$\mathcal{M}$. H.-D. Cao and D.-T. Zhou \cite{Cao-Zhou} further proved that
there exists a positive constant $C$ such that $f$ satisfies the lower
estimate:%
\begin{equation}
f(x)\geq\left(  \frac{1}{2}r(x)-C\right)  ^{2} \label{Sharp lower bound for f}%
\end{equation}
for $x\in\mathcal{M}-B\left(  O,C\right)  $ (see Fang, Man, and Zhang
\cite{FangManZHang} for related estimates).

Define the functions%
\[
\operatorname{V}:\mathbb{R}\rightarrow\lbrack0,\infty),\qquad\operatorname{R}%
:\mathbb{R}\rightarrow\lbrack0,\infty)
\]
by%
\[
\operatorname{V}(c)\doteqdot\int_{\{f<c\}}d\mu,\qquad\operatorname{R}%
(c)\doteqdot\int_{\{f<c\}}R\,d\mu.
\]
In \cite{Cao-Zhou}, the following \textsc{ode} relating $\operatorname{V}(c)$
and $\operatorname{R}(c)$ was established
\begin{equation}
0\leq\frac{n}{2}\operatorname{V}(c)-\operatorname{R}(c)=c\operatorname{V}%
^{\prime}(c)-\operatorname{R}^{\prime}(c). \label{equation number 1}%
\end{equation}

Recall that the asymptotic volume ratio ($\operatorname{AVR}$) of a complete
noncompact Riemannian manifold $(\mathcal{N}^{n},h)$ is defined by%
\begin{equation}
\operatorname{AVR}(h)\doteqdot\lim_{r\rightarrow\infty}\frac
{\operatorname{Vol}B(p,r)}{\omega_{n}r^{n}} \label{AVR display}%
\end{equation}
if the limit exists, where $B(p,r)$ denotes the geodesic ball in $\mathcal{N}$
with center $p$ and radius $r$ and where $\omega_{n}$ is the volume of the
unit Euclidean $n$-ball. It is easy to check that the $\operatorname{AVR}(h)$
is independent of the choice of $p$. Moreover, if $h$ has nonnegative Ricci
curvature, then this limit exists by the Bishop--Gromov volume comparison theorem.

H.-D. Cao and D.-T. Zhou \cite{Cao-Zhou} proved the following using
(\ref{equation number 1}) and aided by an observation of Munteanu \cite{Mu}.

\begin{theorem}
Any complete noncompact shrinking gradient Ricci soliton must have at most
Euclidean volume growth, i.e., $\limsup_{r\rightarrow\infty}\frac
{\operatorname{Vol}B(O,r)}{\omega_{n}r^{n}}$ is finite.
\end{theorem}

Note that an earlier result by Carrillo and Ni \cite{Ca-Ni} states that any
nonflat shrinker with nonnegative Ricci curvature must have zero
$\operatorname{AVR}$. Based on Cao and Zhou's work, Shijin Zhang \cite{Zhang}
proved a sharp upper bound on the volume growth of shrinkers under the
assumption that $R\geq\delta$ for some constant $\delta>0$. More recently,
C.-W. Chen \cite{ChenC} proved that the $\operatorname{AVR}$ of a shrinker is
bounded from below by some $c>0$ if the average scalar curvature satisfies
$\frac{1}{\operatorname{Vol}B(O,r)}\int_{B(O,r)}R\,d\mu\leq r^{\alpha}$, where
$\alpha$ is a negative constant (see also \cite{Cao-Zhou} for a similar result
in the case where $\alpha=0$).

We observe that the results in \cite{Ca-Ni}, \cite{Cao-Zhou}, \cite{Mu},
\cite{Zhang}, and \cite{ChenC} lead to a necessary and sufficient condition
for noncompact shrinkers to have positive $\operatorname{AVR}$.

\begin{theorem}
Let $(\mathcal{M}^{n},g,f)$ be a complete noncompact shrinking gradient Ricci
soliton. Then $\operatorname{AVR}(g)$ exists \emph{(}and is finite\emph{)}.
Moreover, $\operatorname{AVR}(g)>0$ if and only if $\int_{n+2}^{\infty}%
\frac{\operatorname{R}(c)}{c\operatorname{V}(c)}dc<\infty$.
\end{theorem}

\begin{proof}
Let $\operatorname{P}(c)\doteqdot\frac{\operatorname{V}\left(  c\right)
}{c^{\frac{n}{2}}}-\frac{\operatorname{R}\left(  c\right)  }{c^{\frac{n}{2}%
+1}}$ and $\operatorname{N}\left(  c\right)  \doteqdot\frac{\operatorname{R}%
(c)}{c\operatorname{V}(c)}$. Note that $\frac{\operatorname{R}\left(
c\right)  }{\operatorname{V}\left(  c\right)  }$ is the average scalar
curvature over the set $\left\{  f<c\right\}  $. The \textsc{ode}
(\ref{equation number 1}) implies%
\begin{equation}
\operatorname{P}^{\prime}(c)=-\left(  1-\frac{n+2}{2c}\right)  \frac
{\operatorname{R}\left(  c\right)  }{c^{\frac{n}{2}+1}}=-\frac{\left(
1-\frac{n+2}{2c}\right)  \operatorname{N}\left(  c\right)  }%
{1-\operatorname{N}\left(  c\right)  }\operatorname{P}(c).
\label{equation number 3}%
\end{equation}
Since $0\leq\operatorname{R}\left(  c\right)  \leq\frac{n}{2}\operatorname{V}%
(c)$ by (\ref{equation number 1}), we have%
\begin{equation}
\left(  1-\frac{n}{2c}\right)  \frac{\operatorname{V}(c)}{c^{\frac{n}{2}}}%
\leq\operatorname{P}(c)\leq\frac{\operatorname{V}(c)}{c^{\frac{n}{2}}}.
\label{equation number 6}%
\end{equation}
Hence, by the bounds (\ref{Sharp upper bound for f}) and
(\ref{Sharp lower bound for f}) for $f$,%
\[
2^{n}\omega_{n}\operatorname{AVR}(g)=\lim_{c\rightarrow\infty}\frac
{\operatorname{V}(c)}{c^{n/2}}=\lim_{c\rightarrow\infty}\operatorname{P}(c),
\]
which exists by (\ref{equation number 3}).

Integrating (\ref{equation number 3}) yields%
\begin{equation}
\operatorname{P}(c)=\operatorname{P}(n+2)\ e^{-\int_{n+2}^{c}\frac{\left(
1-\frac{n+2}{2c}\right)  \operatorname{N}\left(  c\right)  }%
{1-\operatorname{N}\left(  c\right)  }dc} \label{equation number 4}%
\end{equation}
for $c\geq n+2$. From $\frac{\operatorname{R}(c)}{\operatorname{V}(c)}%
\leq\frac{n}{2}$ it is easy to see that for any $c\in\lbrack n+2,\infty)$ we
have%
\begin{equation}
\frac{1}{2}\int_{n+2}^{c}\operatorname{N}\left(  c\right)  dc\leq\int
_{n+2}^{c}\left(  1-\frac{n+2}{2c}\right)  \frac{\operatorname{N}\left(
c\right)  }{1-\operatorname{N}\left(  c\right)  }dc\leq2\int_{n+2}%
^{c}\operatorname{N}\left(  c\right)  dc. \label{equation number 5}%
\end{equation}
If $\int_{n+2}^{\infty}\operatorname{N}\left(  c\right)  dc=\infty$, then by
(\ref{equation number 4}) we have $\operatorname{AVR}(g)=\frac{1}{2^{n}%
\omega_{n}}\lim_{c\rightarrow\infty}\operatorname{P}(c)=0$.

If $\int_{n+2}^{\infty}\operatorname{N}\left(  c\right)  dc<\infty$, then by
(\ref{equation number 4}) and (\ref{equation number 5}), we have%
\[
\operatorname{P}(c)\geq\operatorname{P}(n+2)\ e^{-2\int_{n+2}^{\infty
}\operatorname{N}\left(  c\right)  dc}>0.
\]
Hence $\operatorname{AVR}(g)>0$.
\end{proof}

\markright{A CONDITION FOR RICCI SHRINKERS TO HAVE POSITIVE AVR}

\textbf{Acknowledgment.} We would like to thank Lei Ni for his interest and encouragement.

\end{document}